\documentclass[12pt]{article}

\usepackage[left=1.25in,right=1.25in,top=1.25in,bottom=1.25in]{geometry}
\usepackage{amsmath,amsfonts,amssymb,amsthm,midpage, mdwlist}
\usepackage{color}

\DeclareMathOperator{\lcm}{lcm}

\DeclareMathOperator{\modulo}{\medspace mod}

\DeclareMathOperator{\sgn}{sgn}

\newcommand{\cm}{\text{,}}
\newcommand{\pd}{\text{.}}

\newcommand{\f}{\mathfrak}
\newcommand{\s}{\mathcal}

\newcommand{\F}{\mathbb{F}}

\newcommand{\N}{\mathbb{N}}

\newcommand{\Q}{\mathbb{Q}}
\newcommand{\R}{\mathbb{R}}

\newcommand{\Z}{\mathbb{Z}}

\newcommand{\rar}{\rightarrow}

\newcommand{\angles}[1]{\langle #1\rangle}
\newcommand{\ceiling}[1]{\lceil #1\rceil}
\newcommand{\floor}[1]{\lfloor #1\rfloor}

\newcommand{\set}[2]{\{\,#1:#2\,\}}

\newtheoremstyle{mystyle}{}{}{\it}{}{\bf}{.}{ }{}
\theoremstyle{mystyle}
\newtheorem{theorem}{Theorem}
\newtheorem{proposition}[theorem]{Proposition}
\newtheorem{definition}[theorem]{Definition}

\newtheorem{corollary}[theorem]{Corollary}

\newtheoremstyle{exstyle}{}{}{}{}{\bf}{.}{ }{}
\theoremstyle{exstyle}

\newtheorem{example}[theorem]{Example}

\begin{document}

\title{Totally real algebraic integers of arboreal height two}
\author{George J. Schaeffer, Stanford University\\{\tt {\small gjschaeff@gmail.com}}}\date{}\maketitle

\begin{abstract}In \cite{Salez}, Salez proved that every totally real algebraic integer is the eigenvalue of some tree. We define the {\it arboreal height} of a totally real algebraic integer $\lambda$ to be the minimal height of a rooted tree having $\lambda$ as an eigenvalue.

In this paper, we prove several results about totally real algebraic integers of arboreal height $\le 2$: We show that all real quadratic integers have arboreal height $\le 2$. We characterize the totally real cubic integers of arboreal height $2$. Finally, we prove that every totally real number field is generated (as a ring over $\Q$) by some $\lambda$ of arboreal height $2$.\end{abstract}

\section{Introduction}

A great deal of research has been done on the connections between the structure of a graph and the {\it geometry} of its adjacency and Laplacian spectra. The most well-known results in {\it spectral graph theory} describe relationships between the sizes of a graph's eigenvalues---or the gaps between them---and its colorability, connectivity, or expansion properties. A concise summary of results in this area is given in \cite{Cai}.

The eigenvalues of a graph\footnote{We mean the eigenvalues of the graph's adjacency matrix, rather than of its Laplacian.} are always {\it totally real}\footnote{An algebraic number is called {\it totally real} if the roots of its minimal polynomial over $\Q$ are all real. For example, $2^{1/3}$ is a real algebraic integer but {\it not totally real}, since the other two roots of $x^3-2$ are imaginary.} algebraic integers. It may therefore be interesting to investigate connections between the structure of a graph and the {\it arithmetic} of its spectrum---the factorization of its characteristic polynomial over $\Q$ and finite fields, the invariants of the number field generated by its eigenvalues (its {\it splitting field}), and so forth. There are a few well-known results in this area, such as a theorem due to Collatz and Sinogowitz that any graph whose characteristic polynomial is irreducible over $\Q$ has a trivial automorphism group \cite[Theorem 10.7]{Yap}. In her recent PhD thesis \cite{Monius}, M\"onius proved several results along these lines, showing for example that every totally real abelian extension of $\Q$ arises as the splitting field of a circulant graph.


In \cite{Estes}, Estes proved a conjecture of Hoffman that every totally real algebraic integer is the eigenvalue of some graph. Salez refined this result by showing that every totally real algebraic integer is the eigenvalue of some tree \cite{Salez}. It is natural to ask which totally real algebraic integers appear in the spectra of graphs when we impose additional structural restrictions.

\begin{definition}Let $\lambda$ be a totally real algebraic integer. The {\em arboreal height} of $\lambda$ is the least $h$ such that there exists a rooted tree of height $h$ (equivalently, an unrooted tree of diameter $2h-1$ or $2h$) having $\lambda$ as an eigenvalue.

We denote the set of all totally real algebraic integers of arboreal height $h$ by $\s{A}_h$.\end{definition}


It is clear that $\s{A}_0=\{0\}$. For arboreal height $1$, we observe that the characteristic polynomial of the {\it star graph} $S_k$ on $k+1$ vertices factors over $\Q$ as $x^{k-1}(x^2-k)$. Thus,
\[\s{A}_1=\set{\pm\sqrt{k}}{k\ge 1}\pd\]

Since trees are bipartite, their characteristic polynomials always have the form $F(x^2)$ or $xF(x^2)$ for some monic $F\in\Z[x]$; consequently, the spectrum of a tree (and each set $\s{A}_h$) is closed under negation. The arboreal height of a given $\lambda$ therefore seems to depend more on $[\Q(\lambda^2):\Q]$ than on $[\Q(\lambda):\Q]$. Indeed, we have already seen that $\lambda\in\s{A}_{1}$ if and only if $\lambda\ne 0$ and $[\Q(\lambda^2):\Q]=1$.

The goal of this paper is to investigate $\s{A}_2$. For convenience, we will let $\angles{k_1^{a_1}\cdots k_n^{a_n}}$ denote the rooted tree obtained by taking a root vertex and connecting it to the central vertices of $a_i$ copies of $S_{k_i}$ for $i=1,\ldots,n$ (with $S_0$ consisting of a single vertex). This tree has $1+\sum_ia_i(k_i+1)$ vertices and its root vertex has degree $\sum_ia_i$.

The set $\s{A}_2$ admits a deceptively simple algebraic description:

\begin{proposition}\label{Height2Criterion}Let $\lambda$ be a totally real algebraic integer with $\lambda^2\notin\Z$. Then, $\lambda\in\s{A}_2$ if and only if there exist nonnegative integers $a_0,a_1,a_2,\ldots$ (almost all zero) satisfying
\[\sum_{k=0}^\infty\frac{a_k}{\lambda^2-k}=1\tag{$*$}\pd\]
When the above holds, $\lambda$ is an eigenvalue of the tree $\angles{\prod_kk^{a_k}}$.\end{proposition}

This description is a consequence of the well-known recursion for the characteristic polynomials of trees \cite[Proposition 5.1.1]{Brouwer}. It is worth noting that the above proposition is once again in terms of $\lambda^2$ rather than $\lambda$. An easy consequence:

\begin{corollary}If $\lambda\in\s{A}_2$ then $d\lambda\in\s{A}_2$ for all nonzero $d\in\Z$. \end{corollary}

\section{Basic tools}

Let $\lambda$ be a totally real algebraic integer with $\lambda^2\notin\Z$ (so that $\lambda\notin\s{A}_0\cup\s{A}_1$), let $n=[\Q(\lambda^2):\Q]$, and let $F(x)$ be the minimal polynomial of $\lambda^2$. Let $\lambda_1,\ldots,\lambda_n\in\R$ be the {\it positive} roots of $F(x^2)$ indexed so that $0<\lambda_1<\cdots<\lambda_n$. Note that for each $i$, $\lambda_i$ is algebraically conjugate to one (or both) of $\pm\lambda$.

\begin{proposition}\label{InterlacingSet}If $\lambda\in\s{A}_2$, then for every $i=1,\ldots,n-1$ there exists a positive $k\in\Z$ satisfying $\lambda_i^2<k<\lambda_{i+1}^2$.\end{proposition}

\begin{proof}Let $T$ be a rooted tree of height $2$ so that $\lambda$ is an eigenvalue of $T$. Since $-\lambda$ is also an eigenvalue of $T$, the real numbers $\lambda_i$ and $\lambda_{i+1}$ must both belong to the spectrum of $T$, and we have $0<\lambda_i<\lambda_{i+1}$.

Let $T'$ be the graph obtained by deleting the root vertex of $T$. By Cauchy's interlacing theorem, there exists an eigenvalue $\mu$ of $T'$ so that $\lambda_i\le\mu\le \lambda_{i+1}$. Since the connected components of $T'$ are all star graphs (or isolated vertices) and $\mu>0$, we must have $\mu\in\s{A}_1$. By our earlier characterization of $\s{A}_1$, it follows that $\mu^2=k$ is a positive integer satisfying $\lambda_i^2\le k\le\lambda_{i+1}^2$. The inequalities are strict because $\lambda^2\notin\Z$.\end{proof}

Let $k_1,\ldots,k_n$ be nonnegative integers with $k_1<\cdots<k_n$. We say that $\{k_1,\ldots,k_n\}$ is a {\it left-interlacing set} for $\lambda^2$ if
\[0\le k_1<\lambda_1^2<k_2<\cdots<\lambda_{n-1}^2<k_n<\lambda_n^2\cm\]
and a {\it right-interlacing set} for $\lambda^2$ if
\[\lambda_1^2<k_1<\lambda_2^2<\cdots<k_{n-1}<\lambda_n^2<k_n\pd\]
Note that $\lambda^2$ has a left-interlacing set if and only if it has a right-interlacing set. Proposition \ref{Height2Criterion} can be rephrased as ``$\lambda$ has arboreal height $2$ only if there is an interlacing set for $\lambda^2$.''

Returning to ($*$) from Proposition \ref{Height2Criterion}, we define the additive monoid
\[\Gamma(\lambda)=\left\{\,\sum_{k=0}^\infty \frac{a_k}{\lambda^2-k}\,:\,\text{$a_0,a_1,a_2,\ldots\in\N$ almost all zero}\,\right\}\pd\]
We have $\lambda\in\s{A}_2$ if and only if $1\in\Gamma(\lambda)$. 

In order to prove that $1\in\Gamma(\lambda)$, it is enough to find two (rational) integers in $\Gamma(\lambda)$ that are relatively prime and of opposite signs. To produce elements of $\Gamma(\lambda)\cap\Z$ we may do the following: Observe that for ($*$) to hold, it must hold upon replacing $\lambda$ with any of $\lambda_1,\ldots,\lambda_n\in\R$. So, choose distinct $k_1,\ldots,k_n\in\N$ and consider the matrix equation
\[\begin{bmatrix}\frac{1}{\lambda_1^2-k_1}&\cdots&\frac{1}{\lambda_1^2-k_n}\\
\vdots&\ddots&\vdots\\
\frac{1}{\lambda_n^2-k_1}&\cdots&\frac{1}{\lambda_n^2-k_n}\\
\end{bmatrix}\begin{bmatrix}v_1\\\vdots\\v_n\end{bmatrix}=\begin{bmatrix}1\\\vdots\\1\end{bmatrix}\pd\]
This equation has a unique solution $\vec{v}(\lambda;k_1,\ldots,k_n)=(v_1,\ldots,v_n)$ with coordinates given explicitly by
\[v_i=-\frac{\prod_{j=1}^n(k_i-\lambda_j^2)}{\prod_{j\ne i}(k_i-k_j)}=-\frac{F(k_i)}{\prod_{j\ne i}(k_i-k_j)}\]
where $F(x)$ is the minimal polynomial of $\lambda^2$.

Note that $\vec{v}(\lambda;k_1,\ldots,k_n)\in \Q^n$ and $\sum_i\frac{v_i}{\lambda^2-k_i}=1$. Thus, if all the $v_i$ are positive, then the least common denominator $\delta(\lambda;k_1,\ldots,k_n)$ of $\{v_1,\ldots,v_n\}$ is a positive element of $\Gamma(\lambda)\cap\Z$. On the other hand, if all the $v_i$ are negative, then $-\delta(\lambda;k_1,\ldots,k_n)$ is a negative element of $\Gamma(\lambda)\cap\Z$.

When $k_1<\cdots<k_n$, we have $\sgn\prod_{j\ne i}(k_i-k_j)=(-1)^{n-i}$, so the denominators in $\vec{v}(\lambda;k_1,\ldots,k_n)$ will alternate signs. In order for all the $v_i$ to have the same sign, $F(k_i)$ with $i=1,\ldots,n$ must alternate signs as well. {\it This is the case precisely when $\{k_1,\ldots,k_n\}$ is an interlacing set for $\lambda^2$.} Summarizing:

\begin{proposition}\label{Summary}Suppose that $\lambda^2$ has an interlacing set $\{k_1,\ldots,k_n\}$, indexed as usual so that $k_1<\cdots<k_n$. Let $\vec{v}(\lambda;k_1,\ldots,k_n)=(v_1,\ldots,v_n)$ be the vector in $\Q^n$ with coordinates
\[v_i=-\frac{F(k_i)}{\prod_{j\ne i}(k_i-k_j)}\]
where $F(x)$ is the minimal polynomial of $\lambda^2$, and let $\delta(\lambda;k_1,\ldots,k_n)$ be the least common denominator of $\{v_1,\ldots,v_n\}$.

Then, $\delta(\lambda;k_1,\ldots,k_n)\in\Gamma(\lambda)\cap\Z$ if $\{k_1,\ldots,k_n\}$ is left-interlacing. On the other hand, $-\delta(\lambda;k_1,\ldots,k_n)\in\Gamma(\lambda)\cap\Z$ if $\{k_1,\ldots,k_n\}$ is right-interlacing. Since $\lambda^2$ has a left-interlacing set if and only if it has a right-interlacing set, $\Gamma(\lambda)\cap\Z$ contains both positive and negative integers.\end{proposition}

It is worth remarking that if $\{k_1,\ldots,k_n\}$ is an interlacing set for $\lambda^2$, then the integer $\lcm\set{\prod_{j\ne i}(k_i-k_j)}{i=1,\ldots,n}$ (or its additive inverse) will be an element of $\Gamma(\lambda)\cap\Z$. However, this fact by itself is not generally helpful: When $n>p$, the pigeonhole principle forces one of the $\prod_{j\ne i}(k_i-k_j)$ to be divisible by $p$. Therefore the LCM will also be divisible by $p$, making it impossible to generate $1$ through the addition of such elements. In short, we also need to consider the {\it numerators} of the coordinates of $\vec{v}(\lambda;k_1,\ldots,k_n)$ when determining whether $\lambda\in\s{A}_2$.

\begin{corollary}\label{Conductor}If there is an interlacing set for $\lambda^2$ then $\Gamma(\lambda)\cap\Z=m\Z$ for some positive integer $m$. We have $\lambda\in\s{A}_2$ if and only if $m=1$.\end{corollary}

\begin{proof}Proposition \ref{Summary} guarantees that if there is an interlacing set for $\lambda^2$, then $\Gamma(\lambda)\cap\Z$ contains both a positive integer and a negative integer. The corollary now follows from the fact that any additive submonoid of $\Z$ that contains a positive integer and a negative integer will be equal to $m\Z$ for some positive integer $m$.\end{proof}


\begin{example}Let $\lambda$ be any root of the octic $f(x)=x^8 - 44x^6 + 567x^4 - 2660x^2 + 3564$. The minimal polynomial of $\lambda^2$ is the quartic $F(x)=f(x^{1/2})$, which is irreducible over $\Q$ with four real roots. Indexing them as described above,
\[0<\underbrace{2.215\ldots}_{\lambda_1^2}<\underbrace{6.813\ldots}_{\lambda_2^2}<\underbrace{9.144\ldots}_{\lambda_3^2}<\underbrace{25.827\ldots}_{\lambda_4^2}\pd\]
So $\{0,3,7,10\}$ is a left-interlacing set for $\lambda^2$ and $\{3,7,10,26\}$ is a right-interlacing set for $\lambda^2$. The corresponding vectors are
\begin{align*}\vec{v}(\lambda;0,3,7,10)&=(\tfrac{594}{35},5,\tfrac{3}{7},\tfrac{8}{5})\text{, and}\\
\vec{v}(\lambda;3,7,10,26)&=(-\tfrac{15}{23},-\tfrac{3}{19},-1,-\tfrac{83}{437})\pd\end{align*}
Giving values of $\delta(\lambda;0,3,7,10)=35=5\cdot 7$ and $\delta(\lambda;3,7,10,26)=437=19\cdot 23$. It follows from Proposition \ref{Summary} that $35,-437\in\Gamma(\lambda)$, and so $1=25\cdot 35+2\cdot(-437)\in\Gamma(\lambda)$.  This proves that $\lambda\in\s{A}_2$.

To construct a rooted tree of height $2$ with $\lambda$ as an eigenvalue, we go back to the vectors and see that
\begin{align*}\tfrac{594}{\lambda^2}+\tfrac{175}{\lambda^2-3}+\tfrac{15}{\lambda^2-7}+\tfrac{56}{\lambda^2-10}&=35\text{, and}\\
\tfrac{285}{\lambda^2-3}+\tfrac{69}{\lambda^2-7}+\tfrac{437}{\lambda^2-10}+\tfrac{83}{\lambda^2-26}&=-437\text{, so}\\
\tfrac{14850}{\lambda^2}+\tfrac{4945}{\lambda^2-3}+\tfrac{513}{\lambda^2-7}+\tfrac{2274}{\lambda^2-10}+\tfrac{166}{\lambda^2-26}&=1\pd\end{align*}
Therefore, $\lambda$ is an eigenvalue of the height-$2$ tree $\angles{0^{14850}\cdot 3^{4945}\cdot 7^{513}\cdot 10^{2274}\cdot 26^{166}}$ on $68231$ vertices.

This is not the smallest tree of height $2$ that has $\lambda$ as an eigenvalue. In fact, $F(x)$ was found during a computer search over small height-$2$ trees as a factor of the characteristic polynomial of $\angles{1^8\cdot 3^4\cdot 8^1\cdot 10^1\cdot 18^3}$, which has $110$ vertices.

\end{example}

%
%

Proposition \ref{InterlacingSet} imposes a sort of ``Archimedean'' restriction on those $\lambda\in\s{A}_2$. There are also ``local'' or mod $p$ restrictions:

\begin{proposition}\label{ModPCondition}If $\lambda\in\s{A}_2$, then for every prime $p$, the reduction of $F(x)$ modulo $p$ must factor into irreducible elements of $\F_p[x]$ of degree $\le p$.\end{proposition}

\begin{proof}Let $\lambda\in\s{A}_2$ and suppose otherwise. Then ($*$) has a solution and there is a prime $p$ so that the reduction of $F(x)$ mod $p$ is divisible by some irreducible $g(x)\in\F_p[x]$ of degree $>p$.

Let $R$ be the ring of integers of $\Q(\lambda^2)$. There exists a prime ideal $\f{p}$ of $R$ above $p$ so that $g(x)$ is the minimal polynomial of $\bar{\lambda}^2$, the image of $\lambda^2$ under the quotient map from the ring $R$ to the residue field $R/\f{p}$.

Since $\deg g(x)>1$, $\bar\lambda^2-k\ne 0$ in $R/\f{p}$ for all $k\in\F_p$. Therefore, each $\frac{1}{\lambda^2-k}$ belongs to the local ring $R_{\f{p}}$, and we may map any solution of ($*$) under $R_{\f{p}}\rar R_{\f{p}}/\f{p}=R/\f{p}$ to obtain
\[\frac{b_0}{\bar\lambda^2}+\cdots+\frac{b_{p-1}}{\bar\lambda^2-(p-1)}=1\]
for some $b_0,\ldots,b_{p-1}\in\F_p$.

Clearing denominators and rearranging, we obtain a polynomial $h(x)\in\F_p[x]$ of degree $p$ such that $h(\bar{\lambda}^2)=0$. But since $\deg h(x)=p<\deg g(x)$, this contradicts the fact that $g(x)$ is the minimal polynomial of $\bar\lambda^2$.\end{proof}

Note that the condition above holds automatically if $p\ge n$, so we only need to check it for the finitely many primes $<n$.

\section{The quadratic case}\label{QuadraticCaseSection}

\begin{theorem}\label{QuadraticCase}If $[\Q(\lambda^2):\Q]=2$, then $\lambda\in\s{A}_2$.\end{theorem}

\begin{proof}We may write $\lambda^2=A+B\sqrt{d}$ where $2A,2B\in\Z$, $B\ne 0$, and $d$ is a positive non-square integer. We have
\[\lambda_2^2-\lambda_1^2=\left|(A+B\sqrt{d})-(A-B\sqrt{d})\right|=2|B|\sqrt{d}>1\cm\]
and therefore there exists a least $k$ satisfying $\lambda_1^2<k<\lambda_2^2$.

The set $\{k-1,k\}$ is left-interlacing for $\lambda^2$. The unreduced denominators of $v_1$ and $v_2$ are $\pm 1$, and therefore $\delta(\lambda;k-1,k)=1\in\Gamma(\lambda)$.\end{proof}

\section{The cubic case}

\begin{theorem}Suppose that $[\Q(\lambda^2):\Q]=3$. We have $\lambda\in\s{A}_2$ iff there is an interlacing set for $\lambda^2$ and $F(x)$ (the minimal polynomial of $\lambda^2$) has a root modulo $2$.
\end{theorem}

\begin{proof}If there is no interlacing set for $\lambda^2$, then $\lambda\notin\s{A}_2$ by Proposition \ref{InterlacingSet}. If the cubic $F(x)$ does not have any roots modulo $2$, then it is irreducible modulo $2$, so $\lambda\notin\s{A}_2$ by Proposition \ref{ModPCondition}.

Suppose then that $\lambda^2$ has an interlacing set and that $F(x)$ has a root modulo $2$. By Corollary \ref{Conductor}, $\Gamma(\lambda)\cap\Z=m\Z$ for some positive integer $m$. Let $k_1=\floor{\lambda_2^2}$ and $k_2=\ceiling{\lambda_2^2}$ so that $k_2=k_1+1$ and
\[\lambda_1^2<k_1<\lambda_2^2<k_2<\lambda_3^2\pd\]

Choose $d\ge 3$ and $k_3>\lambda_3^2$ so that $k_3\equiv k_2+1\equiv k_1+2\modulo d$. Then $\{k_1,k_2,k_3\}$ is a right-interlacing set for $\lambda^2$. The (unreduced) denominators of the coordinates of $\vec{v}(\lambda;k_1,k_2,k_3)$ are
\begin{align*}(k_1-k_2)(k_1-k_3)&\equiv (-1)(-2)\equiv 2\modulo d\cm\\
(k_2-k_1)(k_2-k_3)&\equiv (+1)(-1)\equiv -1\modulo d\text{, and}\\
(k_3-k_1)(k_3-k_2)&\equiv (+2)(+1)\equiv 2\modulo d\pd\end{align*}
The integer $d$ does not divide any of these, so it does not divide $\delta(\lambda;k_1,k_2,k_3)$. Therefore, $d$ does not divide $m$. It follows that $\Gamma(\lambda)\cap\Z$ is either $\Z$ or $2\Z$.

Now we will use the hypothesis that $F(x)$ has a root mod $2$. Since $k_1\not\equiv k_2\modulo 2$, we must have $F(k_1)\equiv 0\modulo 2$ or $F(k_2)\equiv 0\modulo 2$ (or both).

Assume first that $F(k_1)\equiv 0\modulo 2$, so that $F(k_1)$ is even. Choose $k_3'>\lambda_3^2$ so that $k_1-k_3'\equiv 2\modulo 4$. Then $\{k_1,k_2,k_3'\}$ is a right-interlacing set for $\lambda^2$. Consider the coordinates of $\vec{v}(\lambda;k_1,k_2,k_3')$:
\[v_1=-\frac{F(k_1)}{(k_1-k_2)(k_1-k_3')}\pd\]
The numerator is even and the denominator is congruent to $2\modulo 4$, so the reduced fraction will have an odd denominator.
\[v_2=-\frac{F(k_2)}{(k_2-k_1)(k_2-k_3')}\pd\]
The denominator is congruent to $3\modulo 4$, hence odd.
\[v_3=-\frac{F(k_3')}{(k_3'-k_1)(k_3'-k_2)}\pd\]
Since $k_3'\equiv k_1\modulo 2$, $F(k_3')\equiv 0\modulo 2$, so the numerator is even. The denominator is congruent to $2\modulo 4$. Thus, the reduced fraction will have an odd denominator.

In the case where $F(k_2)$ is even, we instead choose $k_3'$ so that $k_2-k_3'\equiv 2\modulo 4$.

We conclude in either case that $\delta(\lambda;k_1,k_2,k_3')$ is odd, so $\Gamma(\lambda)\cap\Z\ne 2\Z$. Therefore, $\Gamma(\lambda)\cap\Z=\Z$ and $\lambda\in\s{A}_2$.\end{proof}

\section{The real cyclotomic case}

An important family of totally real number fields are the {\it real cyclotomic fields} $\Q(\zeta_m+\zeta_m^{-1})$ where $\zeta_m$ is a primitive $m$th root of unity. These are the index-$2$ subfields of the usual cyclotomic fields $\Q(\zeta_m)$.

The {\it real cyclotomic unit} $\zeta_m+\zeta_m^{-1}$ occurs as an eigenvalue of the {\it path graph} on $m-1$ vertices, so its arboreal height is $\le\floor{\frac{m-1}{2}}$. We are interested in the values of $m$ for which $\zeta_m+\zeta_m^{-1}\in\s{A}_2$.

\begin{theorem}\label{RealCyclotomic}The real cyclotomic unit $\zeta_m+\zeta_m^{-1}$ has arboreal height $\le 2$ only if $m\in\{1,2,3,4,5,6,8,10,12,16,20,24\}$.\end{theorem}

\begin{proof}Let $\lambda=\zeta_m+\zeta_m^{-1}$ and let $n=[\Q(\lambda^2):\Q]$. We have $n=1$ for $m\le 4$, and for $m\ge 5$ we have the formula
\[n=\left\{\begin{array}{ll}\frac{1}{4}\phi(m) &\text{if $m$ is divisible by $4$, and}\\
\frac{1}{2}\phi(m)&\text{otherwise.}\end{array}\right.\]
The largest conjugate of $\lambda^2$ is $4\cos(\frac{2\pi}{m})^2\le 4$ with equality holding only if $m\in\{1,2\}$. Therefore, when $n>4$, an interlacing set for $\lambda^2$ does not exist. By Proposition \ref{InterlacingSet} this proves that the arboreal height of $\zeta_m+\zeta_m^{-1}$ is $>2$ in all but finitely many cases:

\begin{itemize}\item We have $n=1$ for $m\in\{1,2,3,4,6,8,12\}$. In these cases, $\lambda^2\in\Z$, so $\lambda\in\s{A}_1$.

\item When $m\in\{5,10,16,20,24\}$ we have $n=2$, so $\lambda\in\s{A}_2$ by Theorem \ref{QuadraticCase}.

\item We have $n=3$ for $m\in\{7,9,14,18,28,36\}$. In these cases, the characteristic polynomial of $\lambda^2$ is cubic and irreducible mod $2$. By Proposition \ref{ModPCondition}, $\lambda\notin\s{A}_2$.

\item We have $n=4$ for $m\in\{15,30,32,40,48,60\}$. When $m\in\{15,30,60\}$, the characteristic polynomial of $\lambda^2$ is irreducible mod $2$, and when $m\in\{32,40\}$, it is irreducible mod $3$. So, $\lambda\notin\s{A}_2$ by Proposition \ref{ModPCondition} for all quartic cases except $m=48$.\end{itemize}

For the remaining case, $\lambda=\zeta_{48}+\zeta_{48}^{-1}$. We will show that $\lambda\notin\s{A}_2$ even though $\lambda^2$ has an interlacing set (e.g., $\{0,1,2,3\}$) and $F(x)=x^4 - 8x^3 + 20x^2 - 16x + 1$ factors as $F(x)\equiv (x+1)^4\modulo 2$ and $F(x)\equiv (x^2+2x+2)^2\modulo 3$. That is, the necessary conditions in Proposition \ref{InterlacingSet} and Proposition \ref{ModPCondition} all hold. Our argument will instead rely on a different mod $3$ condition.

Suppose for contradiction that $\lambda\in\s{A}_2$ so that ($*$) has a solution with $a_0,a_1,a_2,\ldots\in\N$. The set $\{\frac{1}{\lambda^2},\frac{1}{\lambda^2-1},\frac{1}{\lambda^2-2},\frac{1}{\lambda^2-3}\}$ is a basis for $\Q(\lambda^2)$ as a $\Q$-vector space, and the unique solution to $\frac{y_0}{\lambda^2}+\frac{y_1}{\lambda^2-1}+\frac{y_2}{\lambda^2-2}+\frac{y_3}{\lambda^2-3}=1$ has $y_0=1/6$.

For all $k\ge 4$ we can find $b_{k}^{(0)},b_{k}^{(1)},b_k^{(2)},b_k^{(3)}\in\Q$ satisfying
\[\frac{1}{\lambda^2-k}=\frac{b_k^{(0)}}{\lambda^2}+\frac{b_k^{(1)}}{\lambda^2-1}+\frac{b_k^{(2)}}{\lambda^2-2}+\frac{b_k^{(3)}}{\lambda^2-3}\]
By replacing each $\frac{1}{\lambda^2-k}$ for $k\ge 4$ in ($*$) with the above linear combination, we have
\[\frac{a_0+\sum_{k\ge 4}a_kb_k^{(0)}}{\lambda^2}+\frac{a_1+\sum_{k\ge 4}a_kb_k^{(1)}}{\lambda^2-1}+\frac{a_2+\sum_{k\ge 4}a_kb_k^{(2)}}{\lambda^2-2}+\frac{a_0+\sum_{k\ge 4}a_kb_k^{(3)}}{\lambda^2-3}=1\]
Therefore, $a_0+\sum_{k\ge 4}a_kb_k^{(0)}=y_0=1/6$.

We will now prove that the coefficients $b_k^{(0)}$ are actually all $3$-integral. Some routine algebraic manipulation of $F(k)$ shows that for every $k\in\N$ there are unique $c_k^{(0)},c_k^{(1)},c_k^{(2)},c_k^{(3)}\in\Z$ satisfying
\[\frac{1}{\lambda^2-k}=\frac{c_k^{(0)}+c_k^{(1)}\lambda^2+c_k^{(2)}\lambda^4+c_k^{(3)}\lambda^6}{F(k)}\]
and since $F(k)\not\equiv 0\modulo 3$ for all $k\in\N$ (in fact $F(k)\equiv 1\modulo 3$ for all $k\in\N$), each $c_k^{(i)}/F(k)$ is $3$-integral. Expressing each element of $\{\frac{1}{\lambda^2},\frac{1}{\lambda^2-1},\frac{1}{\lambda^2-2},\frac{1}{\lambda^2-3}\}$ in terms of the power basis $\{1,\lambda^2,\lambda^4,\lambda^6\}$ and applying Cramer's rule, we obtain
\[b_k^{(0)}=\frac{\det(M')}{\det(M)}\text{ where }M'=\begin{bmatrix}c_k^{(0)}/F(k)&-3/2&0&-1/2\\c_k^{(1)}/F(k)&13/2&-8&5/2\\c_k^{(2)}/F(k)&-7/2&6&-5/2\\c_k^{(3)}/F(k)&1/2&-1&1/2\end{bmatrix}\text{ and }\det(M)=3\]
Because $\frac{1}{\lambda^2-k}$ is congruent to one of $\frac{1}{\lambda^2-1}$, $\frac{1}{\lambda^2-2}$, $\frac{1}{\lambda^2-3}$ modulo the prime ideal $\f{p}$ above $3$ in $\Q(\lambda^2)$, two columns of $M'$ are equal mod $\f{p}$, so $\det(M')\equiv 0\modulo\f{p}$. Since $\det(M')$ is rational, it follows that $\det(M')\equiv 0\modulo 3$. Therefore, the coefficients $b_k^{(0)}$ for $k\ge 4$ are all $3$-integral. This is a contradiction, since $a_0+\sum_{k\ge 4}a_kb_k^{(0)}=1/6$ with $a_0,a_4,a_5,\ldots\in\N$. \end{proof}

%


\section{General results for number fields}

\begin{theorem}\label{SufficientlyLargeD}Let $\lambda$ be a totally real algebraic integer with $\lambda^2\notin\Z$. There exists a positive rational integer $D$ such that $D\lambda\in\s{A}_2$.\end{theorem}

\begin{proof}First, notice that by replacing $\lambda$ with $d\lambda$ for sufficiently large $d\in\N$, we may assume that there is an interlacing set for $\lambda^2$. Let $\{k_1,\ldots,k_n\}$ be a left-interlacing set for $\lambda^2$ and notice that for any positive integer $D$, we have
\[D^2k_1<(D\lambda_1)^2<D^2k_2<\cdots<D^2k_n<(D\lambda_n)^2\]
so $\{D^2k_1,\ldots,D^2k_n\}$ is left-interlacing for $(D\lambda)^2$.

The minimal polynomial of $(D\lambda)^2$ is $(x-D^2\lambda_1^2)\cdots(x-D^2\lambda_n^2)$ and substituting $x=D^2k$ yields $D^{2n}F(k)$. The $i$th coordinate of $\vec{v}(D\lambda;D^2k_1,\ldots,D^2k_n)$ is given by
\[-\frac{D^{2n}F(k_i)}{\prod_{j\ne i}(D^2k_i-D^2k_j)}=-\frac{D^{2n}F(k_i)}{D^{2n-2}\prod_{j\ne i}(k_i-k_j)}=-D^2\frac{F(k_i)}{\prod_{j\ne i}(k_i-k_j)}\pd\]
Choose $D$ so that $D^2$ cancels each $\prod_{j\ne i}(k_i-k_j)$. Then, $\delta(D\lambda;D^2k_1,\ldots,D^2k_n)=1\in\Gamma(\lambda)$,  so $D\lambda\in\s{A}_2$. \end{proof}

\begin{corollary}Let $K$ be a totally real number field, $K\ne\Q$. There exists $\alpha\in\s{A}_2$ such that $K=\Q(\alpha)$.\end{corollary}

\begin{proof}If $[K:\Q]=2$, then by Theorem \ref{QuadraticCase} we may simply take $\alpha\in K$ with $\alpha^2\notin\Z$.

Otherwise, let $\lambda$ be a primitive element for $K$ so that $K=\Q(\lambda)$. Note that $\lambda^2\notin\Z$. Theorem \ref{SufficientlyLargeD} guarantees that there is a positive integer $D$ large enough so that $\alpha=D\lambda$ has arboreal height $2$, and $K=\Q(\lambda)=\Q(D\lambda)=\Q(\alpha)$.\end{proof}

\begin{corollary}For all $n\ge 2$, there exist $\lambda\in\s{A}_2$ of degree $n$.\end{corollary}

\begin{example}In the proof of Theorem \ref{RealCyclotomic} we saw that $\lambda=\zeta_7+\zeta_7^{-1}=2\cos(\frac{2\pi}{7})$ has arboreal height $3$. The set $\{k_1,k_2,k_3\}=\{0,1,2\}$ is left-interlacing for $\lambda^2$, and the values of $\prod_{j\ne i}(k_i-k_j)$ are $2$, $-1$, and $2$. It follows from the proof of Theorem \ref{SufficientlyLargeD} that $2\lambda$ will have arboreal height $2$.

Indeed, $\{0,4,8\}$ is a left-interlacing set for $4\lambda^2$ with $\vec{v}(2\lambda;0,4,8)=(2,4,2)$. Thus
\[\frac{2}{(2\lambda)^2}+\frac{4}{(2\lambda)^2-4}+\frac{2}{(2\lambda)^2-8}=1\]
so $2\lambda$ is an eigenvalue of the height-$2$ tree $\angles{0^2\cdot 4^4\cdot 8^2}$ on $41$ vertices. Note that while we used Proposition \ref{ModPCondition} to show that $\lambda\notin\s{A}_2$, it does not apply to $2\lambda$, because the minimal polynomial for $4\lambda^2$ is $x^3 - 20x^2 + 96x - 64\equiv x^3\modulo 2$.\end{example}

Every number field is generated by a totally real algebraic integer of arboreal height $2$. Even so, a totally real number field of degree $\ge 3$ must always contain elements of arboreal height $>2$:

\begin{proposition}Let $\lambda$ be a totally real algebraic integer of degree $\ge 3$ so that $\lambda^{-1}$ is also an algebraic integer. At least one of $\lambda$ or $\lambda^{-1}$ has arboreal height $>2$.\end{proposition}

\begin{proof}In this case, the algebraic norm of $\lambda$ is one of $\pm 1$, so we have
\[(\lambda_1^2)\cdots(\lambda_n^2)=(\lambda_1^{-1})^2\cdots(\lambda_n^{-1})^2=1\] Since $n\ge 3$, we have either $\lambda_i^2<1$ and $\lambda_j^2<1$ for distinct $i$ and $j$, or $(\lambda_i^{-1})^2<1$ and $(\lambda_j^{-1})^2<1$ for distinct $i$ and $j$. Thus, either $\lambda$ or $\lambda^{-1}$ cannot have an interlacing set, so either $\lambda$ or $\lambda^{-1}$ has arboreal height $>2$ by Proposition \ref{InterlacingSet}.\end{proof}

\begin{corollary}If $K$ is a totally real number field with $[K:\Q]\ge 3$, then $K$ contains an integer of arboreal height $>2$.\end{corollary}

\begin{proof}This is a consequence of the preceding proposition and Dirichlet's unit theorem for algebraic number fields.\end{proof}

\section{Further cases and limitations}

While it was relatively easy to state necessary and sufficient conditions for $\lambda\in\s{A}_2$ in the cubic case, such conditions become much more intricate for $n\ge 4$. We already saw in the proof of Theorem \ref{RealCyclotomic} that in the quartic case we may have $\lambda\notin\s{A}_2$ even when the necessary conditions given in Propositions \ref{InterlacingSet} and \ref{ModPCondition} all hold.

We conclude with a contrasting example that illustrates another challenge in the quartic case and beyond: When $n\ge 4$, $\Gamma(\lambda)\cap\Z$ is not necessarily generated by elements of the form $\delta(\lambda;k_1,\ldots,k_n)$. That is, we may have $\lambda\in\s{A}_2$ even when the submonoid of $\Gamma(\lambda)\cap\Z$ generated by the various $\delta(\lambda;k_1,\ldots,k_n)$ is not equal to $\Z$.

\begin{example}Let $\lambda$ be any root of the quartic $f(x)=x^4 - x^3 - 24x^2 + 27x + 1$. The minimal polynomial of $\lambda^2$ is $F(x)=x^4-49x^3+632x^2-777x+1$. 

Indexing the roots as usual,
\[0<\underbrace{0.001\ldots}_{\lambda_1^2}<\underbrace{1.370\ldots}_{\lambda_2^2}<\underbrace{23.159\ldots}_{\lambda_3^2}<\underbrace{24.470\ldots}_{\lambda_4^2}\pd\]

Let $\{k_1,k_2,k_3,k_4\}$ be any left-interlacing set for $\lambda^2$ and note that we must have $k_2=1$ and $k_4=24$, meaning that the unreduced denominators of $v_2$ and $v_4$ in $\vec{v}(\lambda;k_0,k_1,k_2,k_3)$ will be multiples of $23$. Since $F(1)=-192$ and $F(24)=-215$, neither of which are divisible by $23$, this factor does not cancel when reducing the fractions $v_0$ and $v_2$, so $\delta(\lambda;k_1,k_2,k_3,k_4)$ will always be a multiple of $23$. A similar statement holds for the right-interlacing sets of $\lambda^2$. Therefore, the submonoid of $\Gamma(\lambda)\cap\Z$ generated by the various $\delta(\lambda;k_1,k_2,k_3,k_4)$ is contained in $23\Z$.

However, we do actually have $\lambda\in\s{A}_2$. We compute
\begin{align*}\vec{v}(\lambda;0,1,2,24)&=(\tfrac{1}{48},\tfrac{192}{\mathbf{23}},\tfrac{599}{44},\tfrac{215}{\mathbf{12144}})\\
\vec{v}(\lambda;0,1,3,24)&=(\tfrac{1}{72},\tfrac{96}{\mathbf{23}},\tfrac{1058}{63},\tfrac{215}{\mathbf{11592}})\\
\vec{v}(\lambda;1,23,24,25)&=(-\tfrac{4}{\mathbf{253}}, -\tfrac{29}{11}, \tfrac{215}{\mathbf{23}}, -12)
\end{align*}
(the boldfaced denominators are divisible by $23$). The key is to notice that we can, in this case, cancel the factors of $23$ in the denominator via linear combination:
\[11\vec{v}(\lambda;0,1,\underline{2},24)+\vec{v}(\lambda;0,1,\underline{3},24)=(\tfrac{35}{144},96,\underline{\phantom{0}}, \tfrac{215}{1008})\]
(we have left the third coordinate on the right blank because the third coordinates of the vectors on the left correspond to coefficients of $\frac{1}{\lambda^2-k}$ for {\it different} values of $k$). So,
\begin{multline*}11\left(\tfrac{1}{48\lambda^2}+\tfrac{192}{23(\lambda^2-1)}+\tfrac{599}{44(\lambda^2-2)}+\tfrac{215}{12144(\lambda^2-24)}\right)\\+\left(\tfrac{1}{72\lambda^2}+\tfrac{96}{23(\lambda^2-1)}+\tfrac{1058}{63(\lambda^2-3)}+\tfrac{215}{11592(\lambda^2-24)}\right)=12\end{multline*}
Collecting like terms and multiplying by $1008=2^4\cdot 3^2\cdot 7$ yields
\[\tfrac{245}{\lambda^2}+\tfrac{96768}{\lambda^2-1}+\tfrac{150948}{\lambda^2-2}+\tfrac{16928}{\lambda^2-3}+\tfrac{215}{\lambda^2-24}=2^6\cdot 3^3\cdot 7\]
hence $12096=2^6\cdot 3^2\cdot 7\in\Gamma(\lambda)$, showing that $\Gamma(\lambda)\cap\Z\not\subseteq 23\Z$.

Now, $\delta(\lambda;1,23,24,25)=253=11\cdot 23$, so $-253$ is also an element of $\Gamma(\lambda)$. Since $12096$ and $-253$ are coprime, $\Gamma(\lambda)\cap\Z=\Z$ and therefore $\lambda\in\s{A}_2$. The height-$2$ tree constructed from the relation $195\cdot12096+9323\cdot(-253)=1$ will have roughly 1.6 billion vertices, and $\lambda$ will be among its eigenvalues.

\end{example}

\end{document}